\documentstyle[10pt]{article}
\newcommand{\Om}{\Omega}

\def\tr{{\rm tr}\,}
\def\Dl{\Delta}
\def\1{^{-1}}
\def\si{\sigma}
\def\al{\alpha}
\def\be{\beta}
\def\La{\Lambda}

\newcommand{\Ibb}[1]{ {\rm I\ifmmode\mkern
            -3.6mu\else\kern -.2em\fi#1}}
\newcommand{\Bbb}[1]{\leavevmode\hbox{\kern.3em\vrule
     height 1.2ex depth -.3ex width .2pt\kern-.3em\rm#1}}

\begin{document}

{
\font\fpgfont=cmr10 scaled 1200
\font\atfont=cmbx10 scaled 1440
\font\atit=cmmi10 scaled 1440
\font\aufont=cmbx10 scaled 1200
\font\abfont=cmr10
\font\abtfont=cmcsc10 scaled 1200

\fpgfont
\begin{center}{\atfont
DEFORMED YANGIANS}
\end{center}
\vskip .5 truecm
\begin{center}
{\atfont AND INTEGRABLE MODELS}
\vskip 1.5 truecm
\end{center}
\begin{center}
{\aufont  P.~P.~Kulish$^a$ and A.~A.~Stolin$^b$}
\end{center}
\vskip 1.5 truecm
\centerline{\hspace{-2mm}$^a$St.Petersburg Department of the Steklov 
Mathematical Institute,}
\smallskip
\centerline{Russian Academy of Sciences,}
\smallskip
\centerline{Fontanka 27, St.Petersburg, 191011, Russia,}
\vskip 1.5 truecm
\centerline{$^b$Department of Mathematics, }
\smallskip
\centerline{Gothenburg University,}
\smallskip
\centerline{S-41296 Goteborg, Sweden}
\smallskip
\bigskip
\centerline{May 1997}

\vskip 1.5 truecm
\centerline{
ABSTRACT}
\medskip
{\advance \parindent by 3mm
\noindent 
Twisted Hopf algebra $sl_\xi(2)$ gives rise to a deformation of the
Yangian ${\cal Y}(sl(2))$. The corresponding deformations of  the integrable 
$XXX$-spin chain and the Gaudin model are discussed.
}
}
\vskip 1.5 truecm

\newpage

\def \arr{\longrightarrow}

\section{Introduction}

The development of the quantum inverse scattering
method (QISM) gave rise to the theory of quantum groups [1--3]: the 
latter influenced the QISM as well (cf. [4]). Taking into account 
recent results on the further deformation of the Yangian algebra $\cal Y$
[5], and applications of the Yangian symmetries to a variety of
integrable models [6], the simplest integrable $XXX$-model 
related to the Yangian ${\cal Y}(gl(2))$ is discussed in this paper.
The deformation related to a simple twist of the Lie algebra
$sl(2):\Delta_\xi={\cal F}\Delta {\cal F}\1$, with an appropriate element
${\cal F}\in U(sl(2))\otimes U(sl(2))$ preserves the regularity
property [7] of the corresponding $R$-matrix [8]: $R(0)={\cal P}$,
where $\cal P$ is the permutation map in $\Bbb C^2\otimes\Bbb C^2$.
This observation enables us to write down a deformed
Hamiltonian for the Heisenberg chain of length $N$ ($XXX_\xi$-model) [8],  
$$
H=\sum\limits_n(\si^x_n\si^x_{n+1}+\si^y_n\si^y_n+
\si^z_n\si^z_{n+1}+\xi^2\si^-_n\si^-_{n+1}+
\xi(\si^-_n-\si^-_{n+1}))\,,
$$
where $\xi$ is a deformation parameter, $\si^x_n$, $\si^y_n$, and 
$\si^z_n$ are Pauli sigma-matrices acting in $\Bbb C^2_n$ related 
to the $n$th site of the chain and
$\si^-_n=\frac12(\si^x_n-i\si^y_n)$. According to the general 
scheme of the QISM one can construct integrable models for other 
values of spin $s= 1, \, \frac3{2}, \, 2, \ldots \,$ as well.

However, it
turns out that this Hamiltonian is not Hermitian, and thus creates 
extra difficulties in constructing the algebraic Bethe
Ansatz $(ABA)$ for this model. Due to the triangular structure
of the twist, we demonstrate that the spectrum of the transfer
matrix and the Bethe equations of the deformed $XXX$-model 
coincide with the usual ones. 

The connection of the $XXX$-spin chain with the representation 
theory of the $sl(2)$ algebra in the framework of the QISM was 
explained in [9]. The deformed $XXX_\xi$-spin chain has a similar 
connection with the representation theory of the twisted $sl(2)$. 

Let us point out that for the periodic boundary condition the obtained 
Hamiltonian coincides with appropriate limit of the seven-vertex model 
Hamiltonian [10]. An interesting application of the Drinfeld twist 
depending on the spectral parameter, to the original $XXX$-model is discussed 
in the recent paper [11]. The Gaudin model directly related to the Yangian,  
was studied also in the framework of the QISM [12]. 

\section{ Yang--Baxter algebra}

Let us consider the Yang--Baxter algebra (see e.g. [7, 13]), 
$$
R(u-v)T_1(u)T_2(v)=T_2(v)T_1(u)R(u-v),
$$
for the entries of the $2\times 2$ matrix
$$
T(u)=\left(
\begin{array}{cc}
 A(u) & B(u) \\
C(u) & D(u)
\end{array}
\right),
$$
given by  twisting the Yang solution,  

\begin{equation}
R(u-v)=F_{21}\left(I-\frac\eta{u-v}\cal P\right) F\1_{12}
=R_\xi-\frac\eta{u-v}\cal P,
\end{equation}

with $\cal P$ as the permutation matrix and 
$$
F_{12}=\left(
\begin{array}{cccc} 
1&0&0&0\\
\xi&1&0&0\\
0&0&1&0\\
0&0&-\xi&1\\
\end{array}
\right),\quad
R_\xi=F_{21}F\1_{12}=
\left(
\begin{array}{cccc}
1&0&0&0\\
-\xi&1&0&0\\
\xi&0&1&0\\
\xi^2&-\xi&\xi&1
\end{array}
\right).
$$
This is a reduction to the fundamental representation $\rho$ in
$\Bbb C^2$ of the universal twist element $\cal F$, $F_{12}= 
(\rho \otimes \rho ){\cal F}$ and 
the universal $R$-matrix $\cal R$ of the twisted Yangian
${\cal Y}_{\eta,\xi}(sl(2))$ [5]. Using the generators $h, \, e$ of the 
$sl(2)$: $[h,\,e]\,=\,-2e$, the universal twist ${\cal F}$ is 
$$
{\cal F}=\, 1+\xi h\otimes e +\frac1{2!} 
\xi^2 h(h+2) \otimes e^2 + \ldots \,.
$$
There are quite a few papers on the non-standard quantum algebra 
$sl_\xi(2)$ (see e.g. [14, 15] and Refs therein). Introducing the variable 
$\sigma : \, 1-2\xi e = \exp(-\sigma)$ this twist is ${\cal F} = 
\exp(h \otimes \si /2 ) $.

To study the corresponding Bethe Ansatz related to 
$T(u)= \{ t_{ij}(u)\}$ , we need all 16
commutation relations $(CR)$ of the form 
$$
\sum t_{ij}(u)t_{kl}(v)=\sum t_{mn}(v)t_{pq}(u)\ldots \,.
$$
Constructing $T(u)=L_N(u)\ldots L_2(u)L_1(u)$ from the local
$L$-operators given by the $R$-matrix (1) itself according to
the QISM [7, 13], we write down at the beginning the CR of $A(u)$
and $D(u)$ with $C(v)$, using notations for the $R$-matrix entries: 
$$
\al(u,v) = 1 + \be(u,v) = 1 - \frac\eta{u-v}\,, 
$$
$$
\begin{array}{l}
A(u)C(v)=\al(u,v)C(v)A(u)-(\be(u,v)C(u)-\xi A(u))A(v)-\\
-\xi D(v)A(u)+(\xi C(v)+\xi^2 D(v))B(u),
\end{array}
$$
$$
\begin{array}{l}
D(u)C(v)=(\al(v,u)C(v)-\xi A(v))D(u)
- (\be(v,u)C(u)-\xi D(u))D(v)+\\
+(\xi C(v)+\xi^2 A(v))B(u),
\end{array}
$$
$$
B(u)C(v)=(C(v)+\xi D(v))B(u)+(\xi B(u)-\be(u,v)D(u))A(v) 
+\be(u,v)D(v)A(u) \, ,
$$
$$
\begin{array}{l}
\al(u,v)C(u)C(v)=(\al(u,v)C(v)-\xi D(v))C(u)+\\
+ (\xi C(v)+\xi^2 D(v))D(u)+(-\xi C(u)-\xi^2 A(u))A(v)+\xi A(u)C(v)\, , 
\end{array}
$$
$$
\al(u,v)A(u)A(v)=(\al(u,v)A(v)-\xi B(v))A(u)+ 
(\xi A(v)+\xi^2 B(v))B(u),
$$
$$
\al(u,v)A(u)B(v)=B(v)A(u)+(\be(u,v)A(v)-\xi B(v))B(u),
$$
$$
B(u)B(v)=B(v)B(u),
$$
$$
A(u)D(v)=D(v)A(u)+(\xi A(u)-\be(u,v)C(u))B(v)\\
+(\be(u,v)C(v)-\xi D(v))B(u),
$$
$$
\al(v,u)D(u)B(v)=B(v)D(u)+(\be(v,u)D(v)-\xi B(v)) B(u), 
$$
$$
\begin{array}{l}
(C(u)+\xi A(u))A(v)=(\al(u,v)A(v)-\xi B(v))C(u)+\\
+(\xi A(v)+\xi^2 B(v))D(u)-\be(u,v)A(u) C(v),
\end{array}
$$
$$
B(u)C(v)=(C(v)+\xi A(v))B(u)+\xi B(u)D(v)+\\
\be(v,u)(A(v)D(u)-A(u)D(v)),
$$
$$
\be(u,v)B(u)C(v)=\be(u,v)B(v)C(u)+(A(v)+\xi B(v))D(u)\\
-(D(u)+\xi B(u))A(v),
$$
$$
D(u)B(v)=\al(u,v)B(v)D(u)-\be(u,v)B(u)D(v)-
B(u)B(v),
$$
$$
\begin{array}{l}
(\al(u,v) D(u) -\xi B(u)) D(v) + (\xi D(u) + \xi^2 B(u))B(v) = \\
= \al(u,v) D(v)D(u)\,.
\end{array}
$$

The first step of the ABA is the assumption of the existence 
of a bare vacuum $\Om$ , such that the monodromy matrix $T(u)$ 
acting on $\Om$, is (lower) triangular
$$
A(u)\,\Om=a(u)\Om,\quad D(u)\,\Om=d(u)\Om,\quad B(u)\,\Om=0.
$$
For the twisted $XXX$-model of spin 1/2 (there are the corresponding states 
for the $XXX_\xi$-models of other spin values too)
\begin{equation}
\Om=\bigotimes\limits^N_{k=1}\left(
\begin{array}{c}
0\\
1
\end{array}
\right),\quad a(u)=1,\quad d(u)=(\al(u))^N.
\end{equation}

Let us study the conditions under which the state vector
$C(v)\Om$ is an eigenvector of the transfer matrix $t(u)$: 
$$
t(u)=\tr T(u)=A(u)+D(u),\quad t(u)\Om=(1+d(u))\Om. 
$$
Using the CR  to intertwine $(A(u)+D(u))$ over $C(v)$ and the
properties of the bare vacuum $\Om$ (2) one gets
\begin{eqnarray*}
t(u)C(v)\Om&=&(1\cdot\al(u,v)+d(u)\al(v,u))C(v)\Om\\
&-& (\be(u,v)+\be(v,u)d(v))C(u)\Om\\
&+&\xi(1-d(u))(1-d(v))\Om \,. 
\end{eqnarray*}
Therefore we can condlude that $C(v)\Om$ is an eigenvector 
of the transfer matrix $t(u)$ if $(\al(v))^N=1\,.$ 
This is the Bethe equation for the one magnon sector. 
The action of $t(u)=A(u)+D(u)$ on the vector $C(v_1)$ $C(v_2)$
$\Om$ results in the following combination of vectors $C\cdot C
\cdot\Om$, $C\Om$ and $\Om$ :
$$
(\al(u,v_1)\al(u,v_2)+d(u)\al(v_1,u)\al(v_2,u))C(v_1)C(v_2)\Om 
+\mbox{``unwanted terms''}\, .
$$
Hence, the eigenvalue $\La(u, \{v_j\}) $ of the transfer matrix $t(u)$ 
is the same as for the $XXX$-model [7, 13], although the structure 
of the ``unwanted terms'' is more complicated and 
the standard Bethe equations, which read 
$$
\left(\frac{v_j-\eta}{v_j}\right)^N=\prod\limits_k'
\frac{v_k-v_j+\eta}
{v_k-v_j-\eta}\,,
$$
are not enough to have
$$
\Psi(v_1,\ldots,v_M)=C(v_1)\ldots C(v_M)\Om
$$
as an eigenvector of the transfer matrix: $t(u) \Psi = \La(u, \{v_j\}) \Psi$. 
This fact is related with triangularity of the 
perturbation of the $XXX$-model Hamiltonian [10]. One can see this 
phenomenon by looking even at the Clebsch-Gordan coefficients of the twisted 
$sl(2)$ algebra in the original basis $h\,, e\,, f$. Due to the connection 
of the twisted coproduct with the standard one by the similarity 
transformation $\Delta_\xi={\cal F}\Delta {\cal F}\1$, the CG coefficients 
of the $sl(2)_{\xi}$ are linear combinations of the usual ones with the 
elements of triangular matrix  $F^{n_1 n_2}_{m_1 m_2} (\pi \,, \rho )$ 
of the twist $ (\pi \otimes \rho ) {\cal F}$ in the tensor product 
of two irreducible representations $V_\pi \otimes V_\rho \, $ (cf.[15]). 

\section{Symmetry algebra $sl_\xi(2)$}

The local structure of the monodromy matrix $T(u)$ leads to
the following asymptotic expansion, 
$$
T_N(u)=L_N(u)\ldots L_1(u)=\prod\limits^N_{k=1} R_{ak}(\xi)
+\frac 1u\sum\limits^N_{k=1} M_k 
{\cal P}_{ak} M_{N-k-1}+O\left(\frac1{u^2}\right), 
$$
where
$M_k=\prod\limits^{k-1}_{m=1}R_{am}(\xi)=(id \otimes \Dl^{(k-1)}) 
R(\xi)$. Using for the constant term the notation 
$$
T_0 = 
\left(
\begin{array}{cc} 
E&0\\
G&E\1
\end{array}
\right)=\prod\limits^N_{k=1}R_{ak}(\xi)=(id \otimes \Dl^{(N-1)}) R(\xi),
$$
one gets the symmetry algebra $sl_\xi(2)$ of the quantum
scattering data\,:  

\begin{eqnarray*}
E\,G&=&GE-\xi\, (1 - E^2),\\
E\,A(u)&=&A(u)E-\xi\, B(u)E,\\
E\,D(u)&=&D(u)E+\xi\, EB(u),\\
E\,B(u)&=&B(u)E,\\
E\,C(u)&=&C(u)E+\xi \,E\,A(u)-\xi\,D(u)E=\\
&=&C(u)E+\xi(A(u)-D(u))E-\xi^2B(u)E;\\ 
G\,B(u)&=&B(u)G-\xi(E\,B(u)+B(u)E\1),\\
G\,A(u)&=&A(u)G-\xi(E\,A(u)-A(u)E\1+B(u)G)+\xi^2B(u)E\1. \\
G\,D(u)&=&D(u)G+\xi(E\,D(u)-D(u)E\1-GB(u))-\xi^2B(u)E. \\
G\,C(u)&=&C(u)G+\xi(E\,C(u)+C(u)E\1-G\,A(u)- D(u)G )+ \\
&+&\xi^2 (D(u)E\1- EA(u)).
\end{eqnarray*}

The element $E$ is the group-like one, and it commutes with the
transfer matrix $t(u)$. The elements $E$ and $G$ have 
the following simple behaviour under the coproduct map: 
$$
\Dl(E)=E\otimes E,\quad \Dl(G)=G\otimes E+E\1\otimes G.
$$
However, to extract from $T(u)$ the third generator of the
twisted Hopf algebra $sl_\xi(2)$  one has to consider the 
$1/u$ term in the asymptotic expansion of $T(u)$ . 

The twist structure of the $R$-matrix (1) leads to the similar 
structure of the monodromy matrix $T(u)$. Hence, the quantum scattering 
data of the $XXX_\xi$-model can be expressed as linear combinations 
of the quantum scattering data of the $XXX$-model with coefficients 
constructed from the twist $\cal F$. 

The twisted $R$-matrix has the same spectral decomposition as
the original one: 
$$
{\cal P} R=F_{12}{\cal P} F\1_{12}=P_+(\xi) - P_-(\xi) \,.
$$
The same is true for the $R$-matrix which depends on the
spectral parameter. Hence, applying the fusion procedure [7] to the
monodromy matrix  one gets relations among the transfer
matrices in different representations of the $sl(2)$ for the auxiliary space 
$$
t^{(l+1)}(u)=t^{(l)}(u)t^{(1)}(u-\eta)+(u - \eta)^N t^{(l-1)}(u). 
$$
Another important ingredient of the QISM is the Baxter
difference equation 
$$
\La(u)Q(u)=Q(u-\eta)+d(u)Q(u+\eta). 
$$
The appropriate quasiclassical limit $\eta \simeq\xi\to0$ gives rise
to a deformation of the Gaudin model, analysed also in the
framework of the QISM [12]. One can apply to the proposed $XXX_\xi$-model 
the functional Bethe Ansatz [16] as well.  

\section*{Acknowledgements}

The authors are grateful to N.~M.~Bogoliubov, I.~V.~Komarov, E.~K.~Sklyanin, 
and V.~Tolstoy for informative discussions. We thank C.Burdik for the 
hospitality in Prague, where this paper was reported at the 6th 
International Colloquium ``Quantum Groups and Integrable Systems''. This 
research was supported by the Swedish National Research Council 
and RFFI grant 96-01-00851.

\end{document}